\documentclass[11pt]{article}
\usepackage{latexsym}
\usepackage{amsmath}
\usepackage{cite}
\usepackage{amsthm,color}
\usepackage{amssymb}
\usepackage{mathrsfs}
\usepackage{times,cases}

\usepackage{lscape}
\newtheorem{theorem}{\noindent Theorem}[section]
\newtheorem{proposition}[theorem]{\noindent Proposition}
\newtheorem{definition}[theorem]{\noindent Definition}
\newtheorem{lemma}[theorem]{\noindent Lemma}
\newtheorem{remark}[theorem]{\noindent Remark}
\newtheorem{corollary}[theorem]{\noindent Corollary}
\numberwithin{figure}{section}
\numberwithin{equation}{section}

\renewcommand{\theequation}{\thesection.\arabic{equation}}

\newcommand{\cA}{{\mathcal A}}

\newcommand{\cE}{{\mathcal E}}

\newcommand{\cG}{{\mathcal G}}
\newcommand{\cH}{{\mathcal H}}

\newcommand{\cM}{{\mathcal M}}

\newcommand{\cW}{{\mathcal W}}

\newcommand{\sA}{{\mathscr A}}

\newcommand{\sH}{{\mathscr H}}
\newcommand{\sK}{{\mathscr K}}
\newcommand{\sM}{{\mathscr M}}

\def\B{\mathbb{B}}
\def\T{\mathbb{T}}

\def\R{\mathbb{R}}

\def\N{\mathbb{N}}
\def\1{\mathbb{1}}
\newcommand{\di}{\mathrm{d}}

\def\bc{\begin{center}}
\def\ec{\end{center}}
\def\be{\begin{equation}}
\def\ee{\end{equation}}
\def\bea{\begin{eqnarray}}
\def\eea{\end{eqnarray}}
\def\ba{\begin{array}}
\def\ea{\end{array}}
\def\benu{\begin{enumerate}}
\def\eenu{\end{enumerate}}
\def\bt{\begin{theorem}}
\def\et{\end{theorem}}
\def\bl{\begin{lemma}}
\def\el{\end{lemma}}
\def\bco{\begin{corollary}}
\def\eco{\end{corollary}}
\def\bn{\begin{numcases}}
\def\en{\end{numcases}}
\def\br{\begin{remark}}
\def\er{\end{remark}}
\def\bd{\begin{definition}}
\def\ed{\end{definition}}
\def\bp{\begin{proposition}}
\def\ep{\end{proposition}}
\def\bo{\begin{proof}}
\def\eo{\end{proof}}
\def\bx{\begin{example}}
\def\ex{\end{example}}
\def\bal{\begin{align}}
\def\eal{\end{align}}

\def\pa{\partial}
\def\al{\alpha}\def\b{\beta}
\def\De{\Delta} \def\de{\delta}
\def\na{\nabla}
\def\lam{\lambda} 

\def\ve{\varepsilon}
\def\Sig{\Sigma}
\def\vp{\varphi}
\def\w{\omega}\def\W{\Omega}
\def\gam{\gamma}\def\Gam{\Gamma}

\def\~{\widetilde}

\def\ol{\overline}

\def\la{\leftarrow}
\def\ra{\rightarrow}

\def\8{\infty}
\def\X{\times}

\def\mb{\mbox}
\def\di{{\rm d}}

\def\es{\emptyset}

\def\sm{\setminus}

\def\ss{\subset}
\def\ssnq{\subsetneqq}

\def\Hs{\hspace{0.8cm}}
\def\hs{\hspace{0.4cm}}

\def\vs{\vskip5pt}

\def\({\left(}
\def\){\right)}

\begin{document}

%%%%%%%%%%%%%%%%%%%%%%%%

\begin{center}
    {\large \bf Recurrent Solutions of a Nonautonomous Modified Swift-Hohenberg Equation}
\vspace{0.5cm}\\
{Jintao Wang$^{1,*}$\quad Lu Yang$^{2}$\quad Jinqiao Duan$^{1,3}$}\\\vspace{0.3cm}

{\small $^{1}$Center for Mathematical Sciences \& School of Mathematics and Statistics,\\ Huazhong
University of Science and Technology, Wuhan, 430074, China\\\vspace{0.2cm}
        $^{2}$School of Mathematics and Statistics, Lanzhou
University,\\ Lanzhou, 730000, China\\\vspace{0.2cm}$^{3}$Department of Applied Mathematics, Illinois Institute of Technology,
Chicago IL 60616, USA}
\end{center}

%%%%%%%%%%%%%%%%%%%%%%%%

\renewcommand{\theequation}{\arabic{section}.\arabic{equation}}
\numberwithin{equation}{section}

%%%%%%%%%%%%%%%%%%%%%%%%%%%%%%%%%%%%%%%%%%%%%%%%%%%%%%%%%%%%%%%%%%%%

%%%%%%%%%%%%%%%%%%%%%%%%%%%%%%%%%%%%%%%%%%%%%%%%%%%%%%%%%%%%%%%%%%%%%%
\begin{abstract}
We consider recurrent solutions of the nonautonomous modified Swift-Hohenberg equation
$$u_t+\De^2u+2\De u+au+b|\na u|^2+u^3=g(t,x).$$
We employ Conley index theory to show that, if the forcing $g:\R\ra L^2(\W)$ is a recurrent function, then there are at least two recurrent solutions in $H_0^2(\W)$ under appropriate assumptions on the parameters $a$, $b$ and $g$.



\end{abstract}

\vspace{-1 cm}

%%%%%%%%%%%%%%%%%%%%%%%%%%%%%%%%%%%%%%%%%%%%%%%%%%%%%%%%%%%%%%%%%%%%
\footnote[0]{\hspace*{-7.4mm}
$^{*}$ Corresponding author.\\
E-mail address: wangjt@hust.edu.cn(J.T. Wang); yanglu@lzu.edu.cn(L. Yang); duan@iit.edu(J.Q. Duan).\\
This work was supported by the NSFC grant (11801190), the
Fundamental Research Funds for the Central Universities Grant
(lzujbky-2018-112). }

%%%%%%%%%%%%%%%%%%%%%%%%%%%%%%%%%%%%%%%%%%%%%%%%%%%%%%%%%%%%%%%%%%%%%%

\section{Introduction}
Periodic and periodic-like phenomena arise in various fields in science and engineering, including the traveling of celestial bodies, fluid systems, migration of animals and recurrence of similar events.
In mathematics, these phenomena are described by time-periodic, quasi-periodic, almost periodic and recurrent motions.
The study of these motions are not only physically important, but also mathematically interesting, especially in the theory of dynamical systems.
Great interests have been attracted by the existence and location of such motions for centuries.

Recurrent motions are a sort of motions more general than periodic and almost periodic ones.
It was first introduced by Birkhoff in \cite{B}, to describe the previously-mentioned ``general'' motion of the so-called ``discontinuous type'' in the phase space in dynamics.
There are important interrelationships between recurrent motions and minimal sets, nonwandering points, Poisson stability for certain $n$-manifold regions (see \cite{B,Se1}).
The problem of existence of recurrent motions for ordinary differential equations was studied by Shcherbakov \cite{Sh2}.

In this paper, we study the existence of recurrent solutions of the following nonautonomous problem,
\bn{}u_t+\De^2u+2\De u+au+b|\na u|^2+u^3=g(t,x),\hs x\in\W,\,t>\tau,\label{1.1}\\
u=\De u=0,\hs x\in \pa\W,\,t\geq\tau,\label{1.2}\\
u(x,\tau)=\varsigma(x),\hs\mb{in }\W,\label{1.3}\en
where $\W$ is an open connected bounded domain in $\R^n$, $a$ and $b$ are arbitrary real constants, $u_t=\frac{\pa u}{\pa t}$, $g$ is the forcing satisfying $g\in L^2_{\rm loc}(\R,L^2(\W))$ and $\varsigma\in L^2(\W)$.
The equation \eqref{1.1} is known in the literature as the modified Swift-Hohenberg equation, and when $b=0$, then \eqref{1.1} is known as the Swift-Hohenberg equation.

The Swift-Hohenberg equation was introduced in 1977 by Swift and Hohenberg (\cite{Swi}) in the research of Rayleigh-B\'enard's convective hydrodynamics (see also \cite{PoM}), arising in geophysical fluid flows in the atmosphere, oceans and the earth's mantle.
It is closely contacted with nonlinear Navier-Stokes equations coupled with the temperature equation.
Later, it has also played a valuable role extensively in the study of plasma confinement in toroidal devices (\cite{LaQ}), viscous film flow, lasers (\cite{LM}) and pattern formation.
In the equation \eqref{1.1}, the modified term $b|\na u|^2$ comes from the study of various pattern formation phenomena involving some kind of phase turbulence or phase transition (\cite{Siva}), which prevents the symmetry $u\ra-u$.

In the previous work, most attention was paid to the existence of attractors (global attractor \cite{Pol,Song}, pullback attractor \cite{Park,Wang}, uniform attractor \cite{XuMa} and random attractor \cite{GuoGuo1,Wang}), bifurcations (dynamical bifurcations \cite{Cho1,Cho2}, nontrivial-solution  bifurcations \cite{Xiao}) and optimal control (\cite{Duan,Sun,Zheng,Zhao}) of different types of modified Swift-Hohenberg equations.
Xiao and Gao in \cite{Xiao} gave specific nontrivial bifurcation solutions that bifurcate from the trivial solution for the modified Swift-Hohenberg equations in rectangular domain in $\R^2$ with periodic boundary value.
Nevertheless, special solutions of \eqref{1.1}, such as (almost) periodic, recurrent solutions, have been sparsely discussed until now.

For recurrent solutions of evolutionary equations, Bongolan-Walsh, Cheban and Duan studied the 2D nonautonomous Navier-Stokes equation with a recurrent external forcing term $g(t)$ in \cite{BCD}
The authors in \cite{BCD} showed that the Navier-Stokes equation has a recurrent solution if $g$ is a recurrent function.
In Li, Wei and Wang's paper \cite{LiW}, the authors defined the locally almost periodicity, which is equivalent to the recurrence for a continuous function used in \cite{BCD}.
In \cite{LiW}, they considered an abstract retarded evolutionary equation with the retarded term $F(u(t-r_1),\cdots,u(t-r_n))$ and the external forcing term $g(t)$, and gave the existence of locally almost periodic solutions when $F$ satisfies the linear growth condition and $g$ is locally almost periodic.
However, to our knowledge, the number (not merely existence) of recurrent solutions has been barely studied yet.
\vs

In this paper, we adopt the recurrent function used in \cite{BCD} (also the locally almost periodic function in \cite{LiW}).
In \cite{BCD}, the authors respectively gave the definitions of recurrence for flows, semiflows and nonautonomous dynamical systems, without expounding the consistence of the different recurrences.
Actually for a flow, if a full solution $\gam$ is recurrent, $\gam$ is also a recurrent function (see Remark \ref{re3.2} below).
Based on this observation, for either semiflows or nonautonomous dynamical systems, by saying that a solution $\gam$ is recurrent, we uniformly mean that $\gam$ is defined all over $\R$ and a recurrent function.

With this realization, we extend the classical Birkhoff Recurrence Theorem to the case of semiflows.
And following this result, we give the main theorem (Theorem \ref{th3.3}) for the existence and location of recurrent solutions of a general nonautonomous differential equation with a recurrent forcing.
It is stated in the main theorem that each compact invariant set of the skew product flow corresponding the nonautonomous differential equation contains a recurrent solution.

We apply the main theorem to the study of recurrent solutions of nonautonomous modified Swift-Hohenberg equations \eqref{1.1}.
Different with the linear growth assumed in \cite{LiW}, the nonlinear term of \eqref{1.1} has a super-linear growth and a modified term related to the gradient of $u$, which increases the difficulty of analysis.
However, the topological tool, Conley index can help to overcome these difficulties.
Thus, we employ theories of Morse decomposition and Conley index to find disjoint nonempty compact invariant sets, and hence obtain different recurrent solutions with distinct initial values by the main theorem.

This paper is organized as follows.
In Section 2, we recall some basic knowledge of dynamical systems, recurrent functions and Conley index theory.
In Section 3, we extend the Birkhoff Recurrence Theorem to semiflows, and give our main theorem concerning the recurrent solutions of nonautonomous differential equations.
In the last section, we consider the $n$-dimensional nonautonomous modified Swift-Hohenberg equation \eqref{1.1} with $n\leq 3$, and prove that there are at least two recurrent solutions under some appropriate assumptions on the parameters $a$, $b$ and the forcing $g$ if $g$ is a recurrent function.

\section{Preliminaries}

We introduce some basic concepts and results of dynamical systems, recurrent functions and Conley index theory.

\subsection{Dynamical systems}
Let $(X,\di)$ be a complete metric space, $\T=\R^+$ or $\R$ and $\Phi$ be an autonomous dynamical system on $X$, i.e., $\Phi:\T\X X\ra X$ is a continuous map such that $\Phi(0,x)=x$ and $\Phi(t+s,x)=\Phi(t,\Phi(s,x))$, for all $x\in X$ and $s,t\in\T$.
If $\T=\R^+$, we call $\Phi$ a {\em semiflow}; if $\T=\R$, we call $\Phi$ a {\em flow}.

A semiflow $\Phi$ is said to be {\em gradient} if there is a continuous function $V:X\ra\R$ such that $t\ra V(\Phi(t,x))$ is non-increasing for each $x\in X$ and if $x$ is such that $V(\Phi(t,x))=V(x)$ for all $t\geq0$, then $x$ is an equilibrium of $\Phi$, i.e., $\Phi(t,x)=x$ for all $t\geq0$.
The function $V$ is called a {\em Lyapunov function} of $\Phi$.

A subset $K\subset X$ is called {\em invariant}, if $\Phi(t)K=K$ for all $t\in\T$.
An invariant set $K$ of $\Phi$ is said to be {\em isolated} if there is a neighborhood $N$ of $K$ such that $K$ is the maximal invariant set in $N$.
Correspondingly, $N$ is called an {\em isolating neighborhood} of $K$.
A subset $K\subset X$ is called {\em minimal}, if it is nonempty, closed and invariant, and it contains no proper subset with these three properties.
A subset $N$ of $X$ is said to be {\em admissible}, if for arbitrary sequences $x_k\in N$ and $t_k\ra+\8$ with $\Phi([0,t_k])x_k\subset N$ for all $k$, the sequence of the end points $\Phi(t_k)x_k$ has a convergent subsequence.

A {\em solution} of $\Phi$ is a map $\gamma:J\ra X$, where $J$ is an interval, such that $\gamma(t)=\Phi(t-s,\gamma(s))$ for all $t,s\in J$ with $t\geq s$.
If $J=\R$, we say $\gam$ is a {\em full} solution.
By the statement ``a solution $\gam:J\ra X$ {\em through} $x$'', we mean $\gam(t)=x$ for some $t\in J$.
If $\gam$ is defined on an interval containing $[0,\8)$, the {\em $\w$-limit set} $\w(\gam)$ is defined as
$$
\w(\gam)=\{y\in X:\,\mb{there exists }t_n\ra\8\mb{ such that }\gam(t_n)\ra y\}.
$$
If $\gam$ is defined on an interval containing $(-\8,0]$, the {\em $\al$-limit set} $\al(\gam)$ is defined as
$$
\al(\gam)=\{y\in X:\,\mb{there exists }t_n\ra-\8\mb{ such that }\gam(t_n)\ra y\}.
$$

Let $\Phi$ be gradient with its Lyapunov function $V$ and $\gam$ be a full solution with $\ol{\gam(\R)}$ compact.
Then we have some simple facts as follows (see \cite{Car}).
\benu\item[(1)] Every compact invariant set contains equilibria, and $\w(\gam)\cup\al(\gam)$ consists of equilibria.
\item[(2)] Both $\al(\gam)$ and $\w(\gam)$ are compact invariant sets and $V(x)\equiv C$ for all $x\in\w(\gam)$ ($x\in\al(\gam)$) and some constant $C$.
\item[(3)] If $\Phi$ has only finitely many equilibria, then $\w(\gam)=\{e\}$ ($\al(\gam)=\{e\}$) for some equilibrium $e$.
\item[(4)] If $e_0$, $\cdots$, $e_k$ are distinct equilibria with $e_0=e_k$, then there are no full solutions $\gam_1$, $\cdots$, $\gam_k$ such that $\al(\gam_{i})=\{e_{i-1}\}$ and $\w(\gam_{i})=\{e_i\}$, for all $i=1$, $\cdots$, $k$.\eenu

A compact invariant set $\cA\ss X$ is said to be an {\em attractor} of $\Phi$, if it {\em attracts} a neighborhood $U$ of itself, that is to say, for all $\ve>0$, there exists $T>0$ such that $\di(\Phi(t,x),\cA)<\ve$, for all $x\in U$ and $t>T$.
The set $B(\cA):=\{x\in X:\,\cA\mb{ attracts }x\}$ is called the {\em attraction basin} of $\cA$.
The attraction basin $B(\cA)$ is always open for each attractor $\cA$.
If $B(\cA)=X$, we say $\cA$ is the {\em global attractor} of $\Phi$ in $X$.

A set $A\ss\cA$ is called an {\em attractor of $\Phi$ in $\cA$}, if it is an attractor of $\Phi|_\cA$, which is the semiflow $\Phi$ restricted on $\cA$.
An ordered collection $\sM=\{M_1,\,\cdots,\,M_k\}$
of subsets $M_k\ss \cA$ is called a {\em Morse decomposition} of $K$, if there exists an increasing sequence
$\es=A_0\ssnq A_1\ssnq\cdots\ssnq A_k=\cA$ of attractors in $\cA$ such that $M_i=A_i\sm B_{\cA}(A_{i-1})$, $1\leq i\leq k$,
where $B_{\cA}$ means the attraction basin for $\Phi|_{\cA}$.
Each $M_i$ is called a {\em Morse set} of $\cA$.

It is known that the Morse set is a compact isolated invariant set (see \cite{Ryb}) and moreover one can give a Morse decomposition in the following way.
\bp[\cite{Ryb}]\label{pro2.1} Let $\sM=\{M_1,\,\cdots,\,M_k\}$ be a family of nonempty, compact, invariant and disjoint subsets of the attractor $\cA$.
Suppose that for each full solution $\gam$ in $\cA$, either $\gam(\R)\ss M_i$ for some $i$ or else there are indices $i<j$ such that $\w(\gam)\ss M_i$ and $\al(\gam)\ss M_j$.
Then $\sM$ is a Morse decomposition of $\cA$.\ep
\subsection{Recurrent functions}

Let $\Phi$ be a flow on $X$.
For each $x\in X$, we can define the {\em hull} of $x$ as
$$\cH(\Phi,x)=\ol{\{\Phi(t,x):t\in\R\}}$$
with the closure taken in $(X,\di)$.
If $\Phi$ is clear, we simply write $\cH(x)=\cH(\Phi,x)$.

A point $x\in X$ is said to be ({\em Birkhoff}) {\em recurrent} for the flow $\Phi$, provided that
\benu\item[(1)] for every $\ve>0$, there is a positive number $l$, such that $\di(x,S([a,a+l])x)<\ve$, for all $a\in\R$; and
\item[(2)] $\cH(x)$ is compact in $X$.\eenu
The corresponding full solution $\Phi(t,x)$ is called a {\em recurrent solution}.
Concerning the minimality and recurrence, we have an simple deduction of Birkhoff Recurrence Theorem in \cite{Se1}.
\bt\label{th2.1} Let $\Phi$ be a flow on $X$ and $A\ss X$.
Then $A$ is a compact minimal set if and only if every $x\in A$ is a recurrent point with $A=\cH(x)$.\et

Denote by $C(\R,X)$ the set of all continuous functions $g:\R\ra X$, equipped with the compact-open topology, which is metrizable with the complete metric $\rho$ such that
$$\rho(g_1,g_2)=\sum_{n=1}^{\8}\frac1{2^n}
\frac{\max_{t\in[-n,n]}\di(g_1(t),g_2(t))}{1+\max_{t\in[-n,n]}\di(g_1(t),g_2(t))},$$
for $g_1,g_2\in C(\R,X)$.

In the sequel we always denote $S:=C(\R,X)$.
Let $\theta$ be the shift on $S$, i.e., $\theta_\tau g=g(\cdot+\tau)$, $g\in S$.
It is easy to know that $\theta$ is continuous on $\R\X S$.
Thus $\theta$ determines a flow on $S$, which is called {\em Bebutov's dynamical system}.
For sake of notational simplification, we will use the same $\theta$ for different spaces $X$, and denote the hull $\cH(\theta,g)$ by $\sH(g)$ for every $g\in S$.

\bd A continuous function $g\in S$ is said to be {\bf recurrent} provided that $g$ is a recurrent point for $\theta$ in $S$.\ed

\br In \cite{LiW}, the authors gave a concept --- locally almost periodicity, for a continuous function.
A continuous function $g:\R\ra X$ is said to be locally almost periodic, if for every $\ve,T>0$, there exists $l>0$, such that for all $a,s\in\R$, there exists $\tau_s\in[a,a+l]$, such that
$$\max_{t\in[-T,\,T]}\di(\theta_{s}g(t),\theta_{\tau_s}g(t))<\ve.$$

This concept is indeed equivalent to recurrence for $g\in S$, by Theorem 16 in \cite{LiW} and the Birkhoff Recurrence Theorem (Theorem \ref{th2.1}).
Locally almost periodicity can be viewed as another description of recurrence, without dynamical systems involved.\er

We have the following consequence according to the properties of locally almost periodicity in Proposition 15 of \cite{LiW}.
\bp\label{pro2.3} A recurrent function $g$ is uniformly continuous on $\R$ such that $\ol{g(\R)}$ is compact in $X$.\ep
\subsection{Conley index}
We now recall the Conley index theory, see \cite{Ryb} for details.

Let $\Phi$ be a semiflow on $X$ and $N$, $E$ be two closed subsets of $X$.
The subset $E$ is said to be {\em $N$-positively invariant}, if for all $x\in E\cap N$ and $t\geq0$, we have $\Phi([0,t))x\ss E$ whenever $\Phi([0,\,t))x\ss N$.
The subset $E$ is said to be an {\em exit set} of $N$, if $E$ is $N$-positively invariant and for every $x\in N$ with some $t<\8$ such that $\Phi(t)x\notin N$, there exists $t_0\leq t$ such that $\Phi(t_0)x\in E$.

Let $K$ be a compact isolated invariant set.
A {\em Conley index pair} $(N,E)$ of $K$ is a closed pair such that $E$ is an exit set of $N$, $\ol{N\setminus E}$ is admissible and $\ol{N\sm E}$ is an isolating neighborhood of $K$.
For each Conley index pair $(N,E)$ of $K$, $(N/E,[E])$ always has the same homotopy type (see \cite{Ryb}).
Here $N/E$ is defined as follows.

If $E\neq\es$, then the space $N/E$ is obtained by collapsing $E$ to a single point $[E]$ in $N\cup E$.
If $E=\es$, we choose a single isolated point $*\notin N$ and define $N/E$ to be the space $N\cup\{*\}$ equipped with the sum topology.
In the latter case we still use the notation  $[E]$ to denote the base point $*$.

\bd The {\bf Conley index} of $K$, denoted by $h(\Phi,K)$, is defined to be the homotopy type of $(N/E,[E])$, where $(N,E)$ is a Conley index pair of $K$.\ed

In this paper, we use $\Sig^r$ to denote the homotopy type of $r$-dimensional pointed sphere and $\ol0$ to denote the homotopy type of pointed one-point set $\{\{*\},*\}$.
A basic fact is that $h(\Phi,K)\ne\ol0$ indicates that $K\ne\es$.

We introduce the continuation property of Conley index in the following.

For a sequence of semiflows $\Phi_k$ on $X$, we write $\Phi_k\ra\Phi_0$, if for all sequences $x_k\in X$ and $t_k\in\R^+$ with $x_k\ra x_0$ and $t_k\ra t_0$, $\Phi_k(t_k)x_k\ra\Phi_0(t_0)x_0$.
A set $N\ss X$ is said to be {\em $\{\Phi_k\}$-admissible} if for two arbitrary sequences $x_k\in X$ and $0\leq t_k\ra\8$ satisfying $\Phi_k([0,t_k])x_k\ss N$ for all $k\in\N^+$, the sequence of endpoints $\Phi_k(t_k)x_k$ has a convergent subsequence.

Let $\cE$ be a metric space. We write $\Phi_\ve\ra\Phi_{\ve_0}$, if $\Phi_{\ve_k}\ra\Phi_{\ve_0}$ for every sequence $\ve_k\in\cE$ with $\ve_k\ra\ve_0$.
The pair $(\Phi_\ve,\,K_\ve)$ is said to be {\em S-continuous at $\ve_0\in\cE$},
if there is a positive $\de$ and a closed subset $N$ of $X$ such that the following two conditions are fulfilled:
\benu\item[(1)] for every $\Phi_\ve$ and $\ve$ with $\rho(\ve,\ve_0)<\de$, the subset $N$ is an admissible closed neighbourhood of $K_\ve$;
\item[(2)] Whenever $\ve_k\ra\ve_0$, then $\Phi_{\ve_k}\ra\Phi_{\ve_0}$ and $N$ is $\{\Phi_{\ve_k}\}$-admissible.\eenu
If $(\Phi_\ve,\,K_\ve)$ is S-continuous at each point $\ve\in\cE$, $(\Phi_\ve,\,K_\ve)$ is said to be {\em S-continuous on $\cE$}.

By the supposition of S-continuity, Conley index possesses the following property, which is called the {\em continuation property}.
\bt[\cite{Ryb}]\label{th4.9}
Let $K_\ve$ be a compact isolated invariant set of $\Phi_{\ve}$ for each $\ve$ lying in a connected component $\cE_0$ of $\cE$.
Suppose that $(\Phi_\ve,\,K_\ve)$ is S-continuous in $\cE_0$.
Then $h(\Phi_{\ve},K_{\ve})$ is constant for $\ve\in\cE_0$.\et

Let $(X,x_0)$ and $(Y,y_0)$ be two pointed spaces.
The {\em wedge sum} $(X,x_0)\vee(Y,y_0)$ and {\em smash product} $(X,x_0)\wedge(Y,y_0)$ are defined, respectively, as follows,
$$
(X,x_0)\vee(Y,y_0)=(\cW,(x_0,y_0)),\hs (X,x_0)\wedge(Y,y_0)=((X\X Y)/\cW,[\cW]),
$$
where $\cW=X\X\{y_0\}\cup\{x_0\}\X Y$.
Denote $[(X,x_0)]$ the homotopy type of a pointed space $(X,x_0)$.
Since the operations ``$\vee$'' and ``$\wedge$'' preserve homotopy equivalence relations, they can be naturally extended to the homotopy types of pointed spaces.

Let $\Phi_i$ be a semiflow on $X_i$, $i=1$, $2$.
We define the product semiflow $\Phi_1\X\Phi_2$ on $X_1\X X_2$ such that
$$(\Phi_1\X\Phi_2)(t,x_1,x_2)=(\Phi_1(t,x_1),\Phi_2(t,x_2)).$$
Let $K_i$ be a compact isolated invariant set of $\Phi_i$.
Then
\be\label{2.2}h(\Phi_1\X\Phi_2,K_1\X K_2)=h(\Phi_1,K_1)\wedge h(\Phi_2,K_2).\ee

\section{Recurrence for semiflows and nonautonomous dynamical systems}

Let $(X,\di)$ be a complete metric space and $\Phi$ be a semiflow on $X$.
A {\em recurrent solution} of $\Phi$ means a full solution which is recurrent as a continuous function.

Now we extend Theorem \ref{th2.1} to the case of the semiflows in the following theorem.
\bt\label{th3.1}A subset $A$ of $X$ is a compact minimal set if and only if for every point $x\in A$, there is a recurrent solution $\gam$ through $x$ such that $A=\ol{\gam(\R)}$.\et
\bo We first show the sufficiency.
Let $\gam$ be a recurrent solution and $A=\ol{\gam(\R)}$.
By Proposition \ref{pro2.3} and Theorem 16 in \cite{LiW}, we know that $A$ is compact and invariant.
Therefore we only need to show the minimality of $A$, i.e., there are no proper compact invariant sets in $A$.

We prove this by contradiction.
We assume that there is a full solution $\gam_0$ in $A$ such that $A_0:=\ol{\gam_0(\R)}$ is a proper compact invariant set in $A$.
Then by the density of $\gam(\R)$ in $A$, we have a $t_0\in\R$ such that
\be\label{3.1}\gam(t_0)\notin A_0.\ee
Hence $\theta_{t_0}\gam\notin\sH(\gam_0)$; otherwise, there is a sequence $t_n\in\R$ such that $\theta_{t_n}\gam_0\ra\theta_{t_0}\gam$ in $(S,\rho)$, which implies that $\gam_0(t_n)\ra\gam(t_0)$, contradicting \eqref{3.1}.
This indicates that $\sH(\gam_0)$ is a nonempty closed invariant proper subset of $\sH(\gam)$ for $\theta$, which contradicts the minimality of $\sH(\gam)$ (by Theorem \ref{th2.1}).
\vs

Now we verify the necessity.
Let $A$ be a compact minimal set of $\Phi$ and $x\in A$.
Pick an arbitrary full solution $\gam$ through $x$ in $A$ and it is obvious that $A=\ol{\gam(\R)}$.
If $\gam$ is a recurrent solution, we are done; if not, we aim to find a recurrent solution through $x$ via $\gam$.

First we prove the compactness of $\sH(\gam)$, for which we only need to show the pre-compactness of $\sH':=\{\theta_t\gam:t\in\R\}$ in $(S,\rho)$.
We claim that
\be\label{3.2} \gam\mb{ is uniformly continuous on }\R.\ee
By confirming \eqref{3.2}, we know that all the functions in $\sH'$ are equi-continuous on $\R$ and the pre-compactness of $\sH'$ is immediately obtained by Ascoli-Arzela theorem.

We show the claim \eqref{3.2} by contradiction and assume there are $\ve_0>0$ and two sequences $s_n,\,t_n\in\R$ with $t_n>s_n$ and $t_n-s_n\ra0$ as $n\ra\8$, such that
\be\label{3.3}\di(\gam(t_n),\gam(s_n))>\ve_0.\ee
Since $\gam$ is contained in a compact set, the sequences can be taken such that $\gam(s_n)\ra x_0$ and $\gam(t_n)\ra y_0$ and by \eqref{3.3}, we obtain that $\di(x_0,y_0)\geq\ve_0$.
By continuity of $\Phi$, we have
$$y_0\la\gam(t_n)=\Phi(t_n-s_n,\gam(s_n))\ra\Phi(0,x_0)=x_0=\lim_{n\ra\8}\gam(s_n),$$
which is a contradiction and ascertain \eqref{3.2}.

Now we have known that $\sH(\gam)$ is compact and invariant for $\theta$.
By the existence theorem of minimal sets (see Theorem II.11 in \cite{Se1}) and Theorem \ref{th2.1}, there is a recurrent function $\gam_0\in S$ such that $\sH(\gam_0)$ is a compact minimal subset of $\sH(\gam)$.
We claim that
\be\label{3.4}\gam_0\mb{ is a full solution of }\Phi\mb{ in }A.\ee
By \eqref{3.4} and the minimality of $A$, we conclude that $A=\ol{\gam_0(\R)}$.
Then there is a sequence $t_n\in\R$ such that $\gam_0(t_n)\ra x$.
Correspondingly, in $\sH(\gam_0)$, we have a subsequence of $t_n$ (still denoted by $t_n$) such that $\theta_{t_n}\gam_0\ra\gam_1\in\sH(\gam_0)$.
Since $\sH(\gam_0)$ is a compact minimal set, by Theorem \eqref{th2.1}, $\gam_1$ is a recurrent function.

Now we show the claim \eqref{3.4}.
Indeed, by the density of $\sH'$ in $\sH(\gam)$, we have a sequence $t_n\in\R$ such that $\theta_{t_n}\gam\ra\gam_0$ in $(S,\rho)$, and hence $\gam(t_n+t)\ra\gam_0(t)$ in $X$ for all $t\in\R$.
Therefore, for all $t,s\in\R$ with $t>s$,
$$\gam_0(t)\la\gam(t_n+t)=\Phi(t-s,\gam(t_n+s))\ra\Phi(t-s,\gam_0(s)),$$
which proves the claim.

Using the claim \eqref{3.4} for $\gam_1$, we also know $\gam_1$ is a full solution in $A$ for $\Phi$.
Note that
$$x\la\gam_0(t_n)=\theta_{t_n}\gam(0)\ra\gam_1(0).$$
Consequently, $\gam_1$ is a recurrent solution through $x$ with $A=\ol{\gam_1(\R)}$.
The proof is complete.
\eo
\br\label{re3.2} If $\Phi$ is a flow, each point $x\in X$ determines a unique full solution in $X$.
Then combining Theorem \ref{th2.1} and \ref{th3.1}, one can easily conclude that, if $\gam$ is a full solution, $\gam$ is a recurrent solution for $\Phi$ if and only if $\gam$ is a recurrent function as a continuous function.\er

Let $X$ be a Banach space.
We consider a general nonautonomous equation as follows,
\be\label{3.5}u_t=f(u)+g(t),\,\, t>\tau;\hs u(\tau)=u_\tau\in X,\ee
where $u$ is the unknown functional, $f:X\supset D(f)\ra X$ and $g:\R\ra X$ are continuous and $D(f)$ is the domain of $f$.
We assume that \eqref{3.5} has a global solution $u(t,\tau;u_\tau)$ for all $u_\tau\in X$, i.e., $u(t,\tau;u_\tau)$ satisfies \eqref{3.5} for all $t\geq\tau$.
A {\em recurrent solution} of \eqref{3.5} means that the solution is defined on $\R$ and recurrent as a continuous function.

We define the solution operator of \eqref{3.5} by a (continuous) {\em cocycle}, $\vp:\R^+\X\sH(g)\X X\ra X,$
i.e., a continuous map satisfying the following conditions,
$$\vp(0,\~g,x)=x,\hs\vp(t+\tau,\~g,x)=\vp(t,\theta_{\tau}\~g,\vp(\tau,\~g,x)),$$
for all $t$, $\tau\in\R^+$, $\~g\in\sH(g)$ and $x\in X$.
According to the classical theory for nonautonomous dynamical systems (see \cite{Car,Kloeden}), we can define a semiflow $\Phi$ on $\sH(g)\X X$ associated with \eqref{3.5} such that
$$\Phi(t,\~g,x):=(\theta_t\~g,\vp(t,\~g,x)),\hs\mb{for all }t\in\R^+ \mb{ and }(\~g,x)\in\W\X X,$$
which is called the {\em skew product flow}.
We always endow $\sH(g)\X X$ with the metric
$$\di_{\sH(g)\X X}((g_1,x_1),(g_2,x_2))=\rho(g_1,g_2)+\di(x_1,x_2)$$
for all $(g_i,x_i)\in\sH(g)\X X$, $i=1,\,2$.
Denote by $\Pi_1$ and $\Pi_2$ the projections from $\sH(g)\X X$ onto $\sH(g)$ and $X$, respectively.
Then we can obtain our {\it main theorem} in this paper for the nonautonomous equation \eqref{3.5}.

\bt[Existence and Location of Recurrent Solutions]\label{th3.3} Suppose that the forcing $g$ is recurrent.
If the skew product flow $\Phi$ has a nonempty compact invariant set $K$, then \eqref{3.5} has a recurrent solution $\gam$ such that
\be\label{3.6}\{(\theta_tg,\gam(t)):t\in\R\}\ss K.\ee
\et
\bo By Theorem II.11 in \cite{Se1}, $K$ contains a compact minimal set $M$ of the system $\Phi$.
Thus $\Pi_1M$ is nonempty, compact and invariant for $\theta$.
By Theorem \ref{th2.1}, the hull $\sH(g)$ is compact minimal and hence $\Pi_1M=\sH(g)$.
Then there exists at least one point $u\in\Pi_2M$ such that $(g,u)\in M$.
Due to Theorem \ref{th3.1}, the skew product flow $\Phi$ has a recurrent solution $\Gam$ in $M$ such that $\Gam(0)=(g,u)$.

Now we show that $\Pi_2\Gam$ is a recurrent function.
Note that $\sH(\Pi_2\Gam)=\Pi_2\sH(\Gam)$ and $\sH(\Pi_2\Gam)$ is compact by Theorem \ref{th2.1}.
Hence we prove it by contradiction and assume that $\sH(\Pi_2\Gam)$ is not minimal.

Indeed, similar to the argument in the proof of Theorem \ref{th3.1}, we can assume that there is a recurrent function $\gam_0\in\sH(\Pi_2\Gam)$ and $t_0\in\R$ such that $\sH(\gam_0)$ is a proper compact minimal subset of $\sH(\Pi_2\Gam)$ and $\theta_{t_0}(\Pi_2\Gam)\notin\sH(\gam_0)$.
The fact that $\gam_0\in\sH(\Pi_2\Gam)$ indicates the existence of a sequence $t_n\in\R$ such that $\theta_{t_n}(\Pi_2\Gam)\ra\gam_0$ in $\sH(\Pi_2\Gam)$.
By compactness of $\sH(\Gam)$, the sequence $t_n$ has a subsequence (still denoted by $t_n$) such that $\theta_{t_n}\Gam\ra\Gam_1:=(\psi,\gam_0)\in\sH(\Gam)$ for some $\psi\in C(\R,\sH(g))$.
However, since $\theta_{t_0}(\Pi_2\Gam)\notin\sH(\gam_0)$, we correspondingly have $\theta_{t_0}\Gam\notin\sH(\Gam_1)$.
This means that $\sH(\Gam_1)$ is a nonempty proper compact invariant subset of $\sH(\Gam)$, contradicting the recurrence of $\Gam$.

Therefore, the solution $\gam:=\Pi_2\Gam$ is a recurrent solution of \eqref{3.5} and obviously satisfies \eqref{3.6}.\eo
\section{Recurrent solutions of the nonautonomous modified Swift-Hohenberg equation}

\subsection{Existence of recurrent solutions}

In this section we consider the nonautonomous modified Swift-Hohenberg equation \eqref{1.1}.

In the following discussion, we always assume that $n\leq 3$ and $g:\R\ra L^2(\W)$ is a recurrent function.
We adopt the theory of semilinear parabolic equations in \cite{Hen} to discuss our issue.

Let $X=L^2(\W)$ with the norm and inner product denoted by $\|\cdot\|$ and $(\cdot,\cdot)$ respectively.
The norm of spaces $L^p(\W)$ is denoted by $\|\cdot\|_{L^p}$ for each $p\geq1$.
Let $L=\De^2$ with its domain
$$D(L)=\{u\in X\cap H^4(\W):u|_{\pa\W}=\De u|_{\pa\W}=0\}.$$
Note that $L$ is a self-adjoint sectorial operator with compact exponent.
Let $\mu_k$ be the eigenvalues of the operator $-\De:X\cap H^2_0\ra X$.
We know
\be\label{4.4}0<\mu_1<\mu_2<\mu_3<\cdots<\mu_k<\cdots\ra\8.\ee
Then the eigenvalues of $L$ are $\mu_k^2$, $k=1$, $2$, $\cdots$.
Let $e^{-Lt}$ be the semigroup generated by $-L$.
We recall the theory of sectorial operators in \cite{Hen} and have that for $\al\in[0,1)$, there are $C_\al>0$ and $\b\in(0,\mu_1^2)$ such that
\be\label{4.5}\|L^{\al}e^{-Lt}L^{-\al}\|\leq C_\al e^{-\b t}\hs\mb{and}\hs\|L^{-\al}e^{-Lt}\|\leq C_\al t^{-\al}e^{-\b t}\hs\mb{for }t>0.\ee
Moreover, we can thus define the functional spaces $X^\al=D(L^\al)$ for $\al\in[0,1)$, whose norm is denoted by $\|\cdot\|_\al$.
Particularly, the functional space $X^{\frac12}=H^2_0(\W)$.

Define $f:X^\al\ra X$ such that $f(u)=2\De u+au+b|\na u|^2+u^3$ and rewrite $g$ as $g(t):\R\ra X$ for each $t\in\R$ with $g(t)(x)=g(t,x)$.
Then the problem \eqref{1.1} -- \eqref{1.3} can rewritten as
\be\label{4.6} u_t+Lu+f(u)=g(t),\,t>\tau,\Hs u(\tau)=\varsigma\in X.\ee
According to the discussions in \cite{Pol,XuMa} and standard methods in \cite{Hen,Se2}, for every $\varsigma\in X^{\frac12}$, the problem \eqref{4.6} possesses a unique, globally defined, mild solution $u$ such that for all $t\geq\tau$,
\be\label{4.7}u(t)=e^{-L(t-\tau)}\varsigma-\int_\tau^te^{-L(t-s)}(f(u(s))-g(s))\di s.\ee

By Proposition \ref{pro2.3} and \eqref{4.4}, we let
$$\cM(g):=\sup_{t\in\R}\|g(t)\|<\8\hs\mb{and}\hs\lam_0:=\min\{\mu_k^2-2\mu_k:k\in\N^+\}.$$
The {\it eventual consequence} of this section is the following theorem for the existence and lower bound of the number of recurrent solutions of the nonautonomous modified Swift-Hohenberg equation \eqref{1.1}.

\bt\label{th4.1} Suppose that $n\leq 3$ and $g$ is a recurrent function.
If $|b|$ is sufficiently small, the problem \eqref{1.1} -- \eqref{1.3} has at least one recurrent solution for a certain initial value.

Furthermore, if in addition, $\cM(g)$ is sufficiently small and $a<-\lam_0$, the problem \eqref{1.1} -- \eqref{1.3} has at least two recurrent solutions with distinct initial values.\et

\subsection{The Proof of Theorem \ref{th4.1}}

In the following analysis, we need the following inequality.

\bl[Gagliardo-Nirenberg inequality \cite{Se2}] Let $\W$ be an open, bounded domain of the Lipschitz class in $\R^n$.
Assume that $1\leq p\leq\8$, $1\leq q\leq\8$, $r\geq1$, $0<\theta\leq1$ and that
\be\label{4.8}k-\frac{n}{p}\leq\theta\(m-\frac{n}{q}\)-(1-\theta)\frac{n}{r}.\ee
Then there is a positive constant $C$ such that one has
$$\|u\|_{W^{k,p}(\W)}\leq C\|u\|_{W^{m,q}(\W)}^{\theta}\|u\|_{L^r(\W)}^{1-\theta},\hs\mb{for all }u\in W^{m,q}(\W).$$
\el

By Theorem \ref{th3.3}, the existence of recurrent solutions depends on that of compact invariant sets for the skew product flow.
We employ Conley index to study this.
For this aim we consider the following auxiliary equation
\be\label{4.9}u_t+Lu+f(u)=\ve\~g(t)\ee
with $\ve\in[0,1]$ and $\~g\in\sH(g)$.

\bl\label{le4.3} There exists a positive number $\~b$ such that, when $|b|\leq \~b$, there is a positive constant $R=R(n,a,\~b,\cM(g))$ such that every full bounded solution $u_\ve$ of \eqref{4.9} in $X^{\frac12}$ satisfies $\|u_\ve(t)\|_{\frac12}<R$ for all $t\in\R$ and $\ve\in[0,1]$.\el
\bo We first give some estimates of the solutions of \eqref{4.9}.
Taking the inner product of \eqref{4.9} with $u$ in $H$ and the inequalities
$$|(2\De u,u)|\leq\frac14\|\De u\|^2+4\|u\|^2,\hs|\ve(\~g(t),u)|\leq\|u\|^2+\frac{\ve^2}4\|\~g(t)\|^2,$$
$$|b(|\na u|^2,u)|=\frac{|b|}{2}|(u^2,\De u)|\leq\frac14\|\De u\|^2+\frac{b^2}4\|u\|_{L^4}^4$$
into consideration, we have
$$\frac{\di}{\di t}\|u\|^2+\|u\|^2+\|\De u\|^2+\frac12(4-b^2)\|u\|_{L^4}^4+(2a-11)\|u\|^2\leq\frac{\ve^2}{2}\|\~g(t)\|^2.$$
Note that $\|u\|^2\leq|\W|\|u\|_{L^4}^2$ and $\|\~g(t)\|\leq \cM(g)$ by the definition of $\sH(g)$,
where $|\W|$ is the Lebesgue measure of $\W$.
When $|b|\leq \~b<2$, we have a positive constant $R_0=R_0(a,\~b,\cM(g))$ such that
\be\label{4.10}\frac{\di}{\di t}\|u\|^2+\|u\|^2+\|\De u\|^2<R_0^2.\ee
Multiplying \eqref{4.10} by $e^t$ and integrating it over $[\tau,t]$ with respect to $t$, we obtain that
\be\label{4.11}\|u(t)\|^2+\int_{\tau}^te^{s-t}\|\De u(s)\|^2\di s<e^{\tau-t}\|u(\tau)\|^2+R_0^2(1-e^{\tau-t}).\ee

Selecting $k=1$, $p=4$, $m=4$, $q=2$, $r=2$, $\theta=\frac{n+4}{16}$ in \eqref{4.8}, we have a certain $C_1>0$ such that
$$\|\na u\|_{L^4}\leq C_1\|\De^2u\|^\frac{n+4}{16}\|u\|^\frac{12-n}{16},$$
and hence by H\"older inequality and Young's inequality, we have $C'_1=C'_1(\~b)>0$ such that
\begin{align}|b(|\na u|^2,\De^2u)|&\leq|b|\|\na u\|_{L^4}^2\|\De^2u\|\notag\\
&\leq C_1|b|\|\De^2u\|^{\frac{n+12}{8}}\|u\|^{\frac{12-n}{8}}\notag\\
&\leq\frac14\|\De^2u\|^2+C'_1\|u\|^{\frac{2(12-n)}{4-n}}.\label{4.12}\end{align}
Selecting $k=0$, $p=6$, $m=4$, $q=2$, $r=2$, $\theta=\frac{n}{12}$ in \eqref{4.8}, we have a certain $C_2>0$ such that
$$\|u\|_{L^6}\leq C_2\|\De^2u\|^\frac{n}{12}\|u\|^\frac{12-n}{12},$$
and again, we have $C'_2>0$ such that
\be\label{4.13}|(u^3,\De^2u)|\leq\|u\|_{L^6}^3\|\De^2u\|\leq C_2\|\De^2u\|^{\frac{n+4}{4}}\|u\|^{\frac{12-n}{4}}\leq\frac14\|\De^2u\|^2+C'_2\|u\|^{\frac{2(12-n)}{4-n}}.\ee
Now we consider the inner product of \eqref{4.9} with $\De^2u$ in $X$.
Combining \eqref{4.12}, \eqref{4.13} and
$$|(2\De u,\De^2u)|\leq\frac14\|\De^2u\|^2+4\|\De u\|^2,\hs|\ve(\~g(t),\De^2u)|\leq\ve^2\|\~g(t)\|^2+\frac14\|\De^2u\|^2,$$
we have
$$\frac{\di}{\di t}\|\De u\|^2+\|\De u\|^2\leq C_3\(\|\De u\|^2+\|u\|^{\frac{2(12-n)}{4-n}}+1\),$$
where $C_3=C_3(a,\~b,\cM(g)):=|9-2a|+2(C'_1+C'_2)+2\cM(g)^2$.
Thus
\be\label{4.14}\|\De u(t)\|^2\leq e^{s-t}\|\De u(s)\|^2+C_3\int_s^te^{\nu-t}\(\|\De u(\nu)\|^2+\|u(\nu)\|^{\frac{2(12-n)}{4-n}}+1\)\di\nu\ee
Integrating \eqref{4.14} over $[\tau,t]$ with respect to $s$ and dividing it by $t-\tau$, we obtain
\be\label{4.15}\|\De u(t)\|^2\leq\(\frac1{t-\tau}+C_3\)\int_\tau^te^{s-t}\|\De u(s)\|^2\di s+C_3\int_\tau^te^{s-t}\(\|u(s)\|^{\frac{2(12-n)}{4-n}}+1\)\di s.\ee
By \eqref{4.11} and \eqref{4.15}, we conclude that when $t-\tau\geq\frac1{C_3}$,
\be\label{4.16}\|\De u(t)\|^2<C_3\(1+2R_0^2+2 e^{\tau-t}\|u(\tau)\|^2+\int_\tau^te^{s-t}\|u(s)\|^{\frac{2(12-n)}{4-n}}\di s\).\ee

We claim that if $u_\ve$ is a bounded full solution of \eqref{4.9} in $X^{\frac12}$ with $\ve\in[0,1]$, we always have $\|u_\ve(t)\|<R_0$ for all $t\in\R$.
Assuming that this claim holds, by \eqref{4.16}, one easily sees that the positive constant $R$ defined as
$$R^2=C_3\(1+4R_0^2+R_0^{\frac{2(12-n)}{4-n}}\)$$
is just what we require in this lemma.

We prove the claim by contradiction.
By the embedding of $X^{\frac12}$ into $X$, $u_\ve$ is also a bounded full solution in $X$.
Hence we assume that there exists $t\in\R$ such that $R_0<\|u_\ve(t)\|\leq R_1:=\sup_{t\in\R}\|u_\ve(t)\|$.
Pick $\tau<t+\ln(\|u_\ve(t)\|^2-R_0^2)-\ln(R_1^2-R_0^2)$.
By \eqref{4.11}, we have
$$\|u_\ve(t)\|^2\leq e^{\tau-t}(\|u_\ve(\tau)\|^2-R_0^2)+R_0^2<\frac{\|u_\ve(\tau)\|^2-R_0^2}{R_1^2-R_0^2}(\|u_\ve(t)\|^2-R_0^2)+R_0^2\leq \|u_\ve(t)\|^2,$$
which causes a contradiction.
The proof is complete.\eo

According to Lemma \ref{le4.3}, we always assume that $|b|\leq\~b$ in the following discussion.

We write the mild solution of \eqref{4.9} in the cocycle form, $u_\ve(t,\~g,\varsigma)$, which satisfies \eqref{4.7} with $\tau$, $g$ therein replaced by $0$, $\ve\~g$ respectively.
We denote by $\B_R$ the open ball centered at the origin with radius $R$ in $X^{\frac12}$.
Let $\Phi_\ve$ be the skew product flow on $\sH(g)\X X^{\frac12}$ corresponding to \eqref{4.9} and denote simply $\Phi=\Phi_1$ for the original problem.

\bl\label{le4.4} Let $\ve_k$ be a sequence in $[0,1]$ with $\ve_k\ra\ve_0\in[0,1]$ as $k\ra\8$.
Then $\Phi_{\ve_k}\ra\Phi_{\ve_0}$ in $X^{\frac12}$ and the product set $\sH(g)\X\ol\B_R$ is $\{\Phi_{\ve_k}\}$-admissible.\el
\bo By replacing $g$ by $\ve_k\~g$, $k\in\N$, in \eqref{4.7}, it is trivial that $\Phi_{\ve_k}\ra\Phi_{\ve_0}$ in $X^{\frac12}$ (see the standard analysis in Theorem I.2.4 in \cite{Ryb} and Theorem 3.4.1 in \cite{Hen}).

Observe that $\sH(g)$ is compact under the metric $\rho$ and \eqref{4.7} is defined all over $[\tau,\8)$.
In order to show the $\{\Phi_{\ve_k}\}$-admissibility and by the compact embedding of $X^\al$ into $X^{\frac12}$ for $\al\in(\frac12,1)$, it is sufficient to prove that whenever $u_{\ve}([0,t],\~g,\varsigma)$ is contained in $\ol\B_R$ for some $(\~g,\varsigma)\in\sH(g)\X X^{\frac12}$ and $t\geq 1$, then $\|u_{\ve}(t,\~g,\varsigma)\|_\al\leq R_\al$ for some positive constant $R_\al$ independent of $\ve$, $\~g$ and $\varsigma$ and some $\al\in(\frac12,1)$.

Indeed, by \eqref{4.9} and \eqref{4.5}, for all $(\~g,\varsigma)\in\sH(g)\X X^{\frac12}$ and each $\al\in(\frac12,1)$,
\begin{align}\|u_{\ve}(t,\~g,\varsigma)\|_\al&\leq\|L^{\al}e^{-Lt}\varsigma\|+\int_0^t\|L^{\al}e^{-L(t-s)}(f(u(s,\~g,\varsigma))-\ve\~g(s))\|\di s\notag\\
&\leq C_\al t^{-\al}e^{-\b t}\|\varsigma\|+C_\al\int_0^ts^{-\al}e^{-\b s}(\|f(u(s,\~g,\varsigma))\|+\cM(g))\di s.\label{4.17}\end{align}
Note also that
\be\|f(u)\|\leq2\|\De u\|+|a|\|u\|+|b|\|\na u\|_{L^4}^2+\|u\|_{L^6}^3.\label{4.18}\ee
By the embedding of $H^2(\W)$ into $L^2(\W)$, $W^{1,4}(\W)$ and $L^6(\W)$ for the dimension $n=1$, $2$, $3$, we can infer from \eqref{4.17} and \eqref{4.18} that there exists $C_4=C_4(n,a,\~b,\cM(g),\al)>0$ such that
$$\|u_{\ve}(t,\~g,\varsigma)\|_\al\leq C_4t^{-\al}e^{-\b t}\|\varsigma\|_{\frac12}+C_4\int_0^ts^{-\al}e^{-\b s}\(1+\|u\|_{\frac12}+\|u\|_{\frac12}^2+\|u\|_{\frac12}^3\)\di s.$$
Hence, when $\|u_{\ve}(s,\~g,\varsigma)\|_{\frac12}\leq R$ for all $s\in[0,t]$ and some $t\geq1$, we have
\be\label{4.19}\|u_{\ve}(t,\~g,\varsigma)\|_\al\leq R_\al:=C_4R+C_4\(1+R+R^2+R^3\)\int_0^\8 s^{-\al}e^{-\b s}\di s.\ee
The proof is finished.\eo

Observe that in the case when $a\geq-\lam_0$, Lemma \ref{le4.3} and \ref{le4.4} guarantee the existence of a nonempty compact invariant set in $\sH(g)\X X^{\frac12}$ for \eqref{4.6}.
By Theorem \ref{th3.3}, we obtain the existence of recurrent solutions for \eqref{1.1} -- \eqref{1.3}.
This proves the first part of Theorem \ref{th4.1}.

However, for the case when $a<-\lam_0$, we can do more detailed analysis and give more than one recurrent solutions for the problem \eqref{1.1} -- \eqref{1.3}.
In the following discussion, we assume $a<-\lam_0$ and adopt Conley index to determine the existence of disjoint compact invariant sets of \eqref{4.6}.

Let $L(a):=\De^2+2\De+a$. We know that the eigenvalues of $L(a)$ are $\lam_k(a):=\mu_k^2-2\mu_k+a$, $k\in\N^+$.
To utilize Conley index theory, we first consider the semiflow $\Psi$ on $X^{\frac12}$ generated by the autonomous equation
\be\label{4.20}u_t+L(a)u+u^3=0.\ee
Applying \eqref{4.11}, \eqref{4.16} and \eqref{4.19} to \eqref{4.20} and the classical existence theorem of global attractors (see \cite{Se2,Tem}, or by the result given in \cite{Pol}), we know that all the bounded full solutions of \eqref{4.20} constitute the global attractor $\cA$ in $X^{\frac12}$ with $\cA\ss\B_R$ by the admissibility (see \cite{Ryb}).
Thus $h(\Psi,\cA)=\Sig^0$.

\bl\label{le4.5} The origin $0$ is an isolated invariant set and
\be\label{4.21}h(\Psi,0)=\Sig^r,\Hs\mb{ where }r=\sum_{\mu_k^2-2\mu_k+a<0}r_k\ee
and $r_k$ is the geometric multiplicity of $\mu_k$.\el

\bo We mainly use Theorem II.3.1 in \cite{Ryb} to calculate the Conley index of $0$.
If $\lam_k(a)\ne0$ for all $k\in\N^+$, \eqref{4.21} follows immediately from Theorem II.3.1 in \cite{Ryb}.
In the following, we consider the case when $\lam_{k_0}(a)=0$ for some $k_0\in\N^+$.

By the statement of Theorem II.3.1 in \cite{Ryb}, we need to consider the reduced equation on the local center manifold of $0$ in $X^{\frac12}$, which can be described as a map $\xi:X_c^{\frac12}\ra(X_c^{\frac12})^{\bot}$,
where $X_c^{\frac12}$ is the eigenspace of $L(a)$ with respect to $0$ and $(X_c^{\frac12})^{\bot}$ is the orthogonal complement of $X_c^{\frac12}$.
Let $E_1$ be the projection from $X^{\frac12}$ to $X_c^{\frac12}$.

By Theorem II.2.1 and II.2.2 in \cite{Ryb}, we are clear that $\xi(u_c)=O(\|u_c\|_{\frac12}^2)$ as $u_c\ra0$ in $X^{\frac12}_c$.
As a result, when $u_c\ra0$ in $X^{\frac12}_c$, the reduced equation of \eqref{4.20} on the center manifold $\{u_c+\xi(u_c):u_c\in X_c^{\frac12}\}$ can be written as
\be\label{4.22}\frac{\di}{\di t}u_c+E_1(u_c+\xi(u_c))^3=\frac{\di}{\di t}u_c+E_1u^3_c+O(\|u_c\|_{\frac12}^4)=0.\ee
Consider the inner product of \eqref{4.22} with $u_c$ in $X$.
We obtain that
$$\frac12\frac{\di}{\di t}\|u_c\|^2+\|u_c\|_{L^4}^4+O(\|u_c\|_{\frac12}^5)=0,$$
in that $(E_1u_c^3,u_c)=(u_c^3,E_1u_c)=(u_c^3,u_c)=\|u_c\|_{L^4}^4$.
Since $X^{\frac12}_c$ is finite-dimensional, we have
$$\frac12\frac{\di}{\di t}\|u_c\|^2\leq-|\W|^{-1}\|u_c\|^4+O(\|u_c\|^5),$$
and hence $0$ is an attractor for the reduced semiflow on the center manifold.
This indicates that $\{0\}$ is an isolated invariant set of $\Psi$ and $h(\Psi,0)=\Sig^0\wedge\Sig^r=\Sig^r$ by Theorem II.3.1 in \cite{Ryb}.
\eo

\bl\label{le4.6} The attractor $\cA$ has a Morse decomposition $\{K_0,\{0\}\}$ with $h(\Psi,K_0)\ne\ol0$.\el
\bo We make use of properties of gradient semiflows to give the proof.

Let $V:X^{\frac12}\ra\R$ be such that
$$V(u)=\frac12(L(a)u,u)+\frac14\int_{\W}u^4\di x.$$
By \eqref{4.20} one sees that
$$\frac{\di}{\di t}V(u(t))=V'(u)\frac{\di}{\di t}u(t)=\(L(a)u+u^3,\frac{\di}{\di t}u(t)\)=-\|L(a)u+u^3\|^2\leq0.$$
If $V(\Psi(t,u))=V(u)$ for all $t\geq0$, we have that $\|L(a)u+u^3\|=0$ and hence $u$ is an equilibrium of $\Psi$.
This implies that $\Psi$ is a gradient semiflow.

We give some necessary estimates of $V$.
Because $\lam_k(a)\geq\lam_0+a$, by the embedding of $L^4(\W)$ into $L^2(\W)$, we have
\be\label{4.23}V(u)\geq\int_\W\(\frac{\lam_0+a}2u^2+\frac14u^4\)\di x=\frac14\int_\W(u^2+\lam_0+a)^2\di x-\frac{|\W|}4(\lam_0+a)^2.\ee
Particularly, since $\lam_0$ can be reached by some $\mu_k$, the equality of \eqref{4.23} can also be reached by the choosing $u$ to be the eigenfunction of $\mu_k$.
Note that $\cA$ is the global attractor, compact.
The function $V$ reaches its minimum on $X^{\frac12}$ at some point $u_0\in\cA$.
By the theory of variational calculus, the point $u_0$ is an equilibrium of \eqref{4.20}.
If $u\in X^{\frac12}$ is an equilibrium, we know $L(a)u+u^3=0$ and then
\be\label{4.24}V(u)=\frac12(L(a)u,u)+\frac14\int_{\W}u^4\di x=-\frac12(u^3,u)+\frac14\int_{\W}u^4\di x=-\frac14\int_{\W}u^4\di x\leq0,\ee
which indicates that $V$ reaches its maximum at $0$ on $\cA$.

When $a<-\lam_0$, by \eqref{4.23}, the minimum of $V$ is negative and hence the equilibrium $u_0\ne0$ by \eqref{4.24}.
Pick a positive number $\de$ such that $\B_\de$ and $\B_{2\de}$ are isolating neighborhoods of $0$.
Then surely $\ol\B_R\sm\B_\de$ is a closed isolating neighborhood.
By the admissibility of $\ol\B_R$, we know that $\Phi$ possesses a compact isolated invariant set $K_0$ in $\ol\B_R\sm\B_\de$.
Then on $K_0$, $V$ can reach a maximum at some point $u^*$ and $u^*$ is an equilibrium; otherwise, $K_0$ contains a full solution $\gam$ through $u^*$ such that
$$\max_{u\in K_0}V(u)\geq V(\al(\gam))>V(u^*),$$
which is a contradiction.
Obviously, $K_0$ contains all the equilibria of $\Phi$ (including $u_0$) except $0$ and $0>V(u^*)\geq V(u_0)$.
Based on this and the continuity of $V$, we can take $0<\de'<\de$ such that $V(u')>V(u^*)$ for all $u'\in\B_{\de'}$.
And $\ol\B_R\sm\B_{\de'}$ is still an isolating neighborhood of $K_0$.

We show that $\{K_0,\{0\}\}$ is a Morse decomposition of $\cA$.
Let $u\in\cA\sm(K_0\cup\{0\})$ and $\gam$ be a full solution through $u$.
By Proposition \ref{pro2.1}, it is sufficient to prove $\al(\gam)=\{0\}$ and $\w(\gam)\ss K_0$.

By the property of gradient semilfows and \eqref{4.24}, $V(\w(\gam))<V(\al(\gam))\leq0$.
This means $0\notin\w(\gam)$ and therefore $\w(\gam)\ss K_0$.
We show $\al(\gam)=\{0\}$ in the rest.

Indeed, if $V(\al(\gam))=0$, we are done; otherwise, since $\al(\gam)$ consists of equilibria, we can deduce that $\al(\gam)\ss K_0$.
Observe that for all $u'\in\B_{\de'}$ and $t\in\R$,
$$V(u')>V(u^*)\geq V(\al(\gam))>V(\gam(t)).$$
We know that $\gam$ is excluded from $\B_{\de'}$.
Note also $\gam$ is contained in $\cA$, which implies that $\gam$ is a full solution in $\ol\B_R\sm\B_{\de'}$.
It can be derived that $\gam(\R)\ss K_0$, resulting in a contradiction with the choice of $\gam$.

Now that $\{K_0,\{0\}\}$ is a Morse decomposition of $\cA$, by definition $K_0$ is an attractor in $\cA$.
By Corollary 5.11 in \cite{Car}, $K_0$ is also an attractor of $\Psi$ on $X^{\frac12}$.
Surely $K_0$ has a closed isolated neighborhood $N_0$ such that $\Psi(t,N_0)\ss N_0$ for all $t\geq0$.
Thus $(N_0,\es)$ is a Conley index of $K_0$ and $h(\Psi,K_0)\ne\ol0$.
The proof is complete.\eo

Now we add the modified term $b|\na u|^2$ to the equation \eqref{4.20}, written as
\be\label{4.25}u_t+L(a)u+b|\na u|^2+u^3=0,\ee
and denote by $\Psi_b$ the semiflow generated by \eqref{4.25}.
By Lemma \ref{le4.3} and \ref{le4.4}, we know when $|b|\leq \~b$, $\Psi_b$ possesses a global attractor $\cA^b$ in $X^{\frac12}$.
For simplicity, we denote $\Xi:=h(\Psi,K_0)$.

\bl\label{le4.7} There is a positive constant $\hat b$ such that if $|b|<\hat b$, $\cA^b$ contains two disjoint compact isolated invariant sets $K_0^b$ and $K_1^b$ such that
\be\label{4.26} h(\Psi_b,K_0^b)=\Xi\hs\mb{and}\hs h(\Psi_b,K_1^b)=\Sig^r.\ee\el
\bo We first show that there is $\hat b_1\in(0,\~b)$ such that whenever $|b|\leq\hat b_1$, $\Psi_b$ has a compact isolated invariant set $K_1^b$ around $0$ and $h(\Psi_b,K_1^b)=\Sig^r$ with $r$ given by \eqref{4.21}.
Given each convergent sequence $b_k\in[\,-\~b,\~b\,]$ with $b_k\ra b_0$, we know that $\Psi_{b_k}\ra\Psi_{b_0}$ by Theorem I.2.4 in \cite{Ryb} and $\B_\de$ is $\{\Psi_{b_k}\}$-admissible with $\de$ given in the proof of Lemma \ref{le4.6}.

We claim that there is $\hat b_1\in(0,\~b)$ such that $\ol\B_\de$ is an isolating neighborhood for $\Psi_b$ with  each $b$ satisfying $|b|\leq\hat b_1$.
Admitting the claim and recalling the continuation property of Conley index, we obtain that the maximal invariant set in $\ol\B_\de$ can be regarded as the $K_1^b$ that we desire above for all $b$ with $|b|\leq\hat b_1$.

The claim can be proved by contradiction with a standard process.
We suppose that there is a sequence $b_k\in(-\hat b_1,\hat b_1)$ with $b_k\ra0$ as $k\ra\8$, which allows a full solution $\gam_k$ in $\ol\B_\de$ for each $\Psi_{b_k}$ such that $\gam_k(0)\in\pa\B_\de$.
Denoting $v_k=\gam_k(-k)$, we have $\Psi_{b_k}([0,k])v_k\ss\ol\B_\de$.
By the $\{\Psi_{b_k}\}$-admissibility of $\ol\B_\de$, we can assume that
$$\gam_k(0)=\Psi_{b_k}(k)v_k\ra v_0\in\pa \B_\de.$$
It can be referred from Theorem I.4.5 in \cite{Ryb} that through $v_0$ there passes a full solution $\gam$ of $\Psi_0$ in $\ol\B_\de$, which contradicts the fact that $\ol\B_\de$ is an isolating neighborhood of $0$ for $\Psi_0$.
This confirms the claim.

With the similar argument for the compact isolated invariant set $K_0$ of $\Psi$ in $\ol\B_R\sm\B_\de$, we also has a $\hat b_0$ such that $\Psi_b$ possesses a compact isolated invariant set $K_0^b$ in $\ol\B_R\sm\B_\de$ for all $|b|\leq\hat b_0$.
Moreover, $h(\Psi_b,K_0^b)=\Xi$.

Finally, we take $\hat b=\min\{\hat b_0,\hat b_1\}$, which satisfies all the requirements of this lemma.
We finish the proof now.
\eo

Now we go on to consider the equation \eqref{4.9}.
Denote by $\Pi_1$ and $\Pi_2$ the projections from $\sH(g)\X X^{\frac12}$ to $\sH(g)$ and $X^{\frac12}$, respectively.

Let $\cH^*$ be the homotopy type of $(\sH(g)\cup\{*\},*)$ with $*\notin\sH(g)$ and sum topology.
Define $\sA_\ve$ to be the maximal compact invariant set of $\Phi_\ve$ in $\sH(g)\X\ol\B_R$ and denote again $\sA=\sA_1$.
By Lemma \ref{le4.3}, we know that $\sA_\ve\ss\sH(g)\X\B_R$ for each $\ve\in[0,1]$.

In the following discussion, we assume $|b|\leq\hat b$ such that we can use the result of Lemma \ref{le4.7}.

\bl\label{le4.8} There is $\hat M>0$ such that if $\cM(g)\leq\hat M$, then $\sA$ contains two disjoint compact isolated invariant sets $\sK_0$ and $\sK_1$ such that
\be\label{4.27} h(\Phi,\sK_0)=\cH^*\wedge\Xi\hs\mb{and}\hs h(\Phi,\sK_1)=\cH^*\wedge\Sig^r.\ee\el

\bo By Lemma \ref{le4.7}, we know $\sA_0=\sH(g)\X\cA_b$ and
\be\label{4.28}h(\Phi_0,\sH(g)\X K_0^b)=\cH^*\wedge\Xi\hs\mb{and}\hs h(\Phi_0,\sH(g)\X K^b_1)=\cH^*\wedge\Sig^r.\ee

Observe in Lemma \ref{le4.4} that $\sH(g)\X\B_R$ is $\{\Phi_{\ve_k}\}$-admissible and $\Phi_{\ve_k}\ra\Phi_{\ve_0}$, for each convergent sequence $\ve_k\in[0,1]$ with $\ve_k\ra\ve_0$.
With a similar argument to that in Lemma \ref{le4.7}, we easily obtain the existence of $\hat\ve$ in $[0,1]$ such that whenever $\ve\leq\hat \ve$, $\sA_\ve$ has two disjoint compact isolated invariant sets $\sK_{0,\ve}$, $\sK_{1,\ve}$ satisfying \eqref{4.27} with $\sK_0$ and $\sK_1$ therein replaced by $\sK_{0,\ve}$ and $\sK_{1,\ve}$, respectively.

Note that all the discussions above rely on the supremum $\cM(g)$, but independent of the nonautonomous forcing $g$ itself.
Hence we can fix some positive constant $M$.
Accordingly we obtain the $\hat\ve=\hat\ve(M)$ such that for all recurrent functions $\mathcal{G}$ with $\cM(\mathcal{G})=M$, whenever $\~\ve\leq\hat \ve$, $\Phi^{\mathcal{G}}_{\~\ve}$ has two disjoint compact isolated invariant sets $\sK_{0,\~\ve}^{\mathcal{G}}$ and $\sK_{1,\~\ve}^{\mathcal{G}}$ in $\sH(\mathcal{G})\X\B_{R}$ such that
\be\label{4.29}h(\Phi_{\~\ve}^{\mathcal{G}},\sK_{0,\~\ve}^{\mathcal{G}})=\cH_{\mathcal{G}}^*\wedge\Xi\hs\mb{ and }\hs h(\Phi_{\~\ve}^{\mathcal{G}},\sK_{1,\~\ve}^{\mathcal{G}})=\cH_{\mathcal{G}}^*\wedge\Sig^r.\ee
Here $\Phi^{\mathcal{G}}_{\~\ve}$ is the skew product flow corresponding to \eqref{4.9} with $\ve\~g$ replaced by $\~\ve\mathcal{G}$ and $\cH_{\mathcal{G}}^*$ is the homotopy type of $\sH(\mathcal{G})\cup\{*\}$.

Now we let $\hat M:=\hat\ve M$ and consider the forcing $g$ with $\cM(g)\leq\hat M$.
We will show $\sA_\ve$ has two disjoint compact isolated invariant sets $\sK_{0,\ve}$, $\sK_{1,\ve}$ satisfying \eqref{4.27} for each $\ve\in[0,1]$ with $\sK_0$ and $\sK_1$ therein replaced by $\sK_{0,\ve}$ and $\sK_{1,\ve}$, respectively.

If $\cM(g)=0$, we have $g\equiv0$ and the conclusion is obvious.
If $\cM(g)>0$, let $\mathcal{G}=\frac{M}{\cM(g)}g$.
We thus have $\cM(\mathcal{G})=M$ and therefore, when $\~\ve\leq\hat\ve$, $\Phi^{\mathcal{G}}_{\~\ve}$ has two disjoint compact isolated invariant sets $\sK_{0,\~\ve}^{\mathcal{G}}$ and $\sK_{1,\~\ve}^{\mathcal{G}}$ in $\sH(\mathcal{G})\X\B_{R}$ satisfying \eqref{4.29}.
Note that $\~\ve\mathcal{G}=\frac{\~\ve M}{\cM(g)}g$.
Let $T:\sH(g)\X X^{\frac12}\ra\sH(\mathcal{G})\X X^{\frac12}$ such that $T(\~g,u)=\(\frac{M}{\cM(g)}\~g,u\)$.
We can deduce that for all $t\geq0$ and $(\~g,u)\in\sH(g)\X X^{\frac12}$,
$$T\circ\Phi_{\frac{\~\ve M}{\cM(g)}}(t,\~g,u)=\Phi_{\~\ve}^{\mathcal{G}}\(t,T(\~g,u)\)\hs\mb{and}\hs\frac{\~\ve M}{\cM(g)}\in\left[0,\frac{\hat\ve M}{\cM(g)}\right]=\left[0,\frac{\hat M}{\cM(g)}\right]\supset[0,1],$$
from $\~\ve\in[0,\hat\ve]$.
Since $T$ is a homeomorphism, by Proposition II.3.2 in \cite{Ryb}, we have that $\Phi_\ve$ has two compact isolated invariant sets $\sK_{0,\ve}$ and $\sK_{1,\ve}$ in $\sA_\ve$ for each $\ve\in[0,1]$ such that
\be\label{4.30}h(\Phi_\ve,\sK_{i,\ve})=h\(\Phi_{\frac{\cM(g)\ve}{M}}^{\cG},\sK_{i,\frac{\cM(g)\ve}{M}}^{\cG}\),\hs i=0,\,1,\ee
due to the fact that $\frac{\cM(g)\ve}{M}\leq\hat\ve$.
Since $\mathcal{G}$ is a positive multiple of $g$, one can easily infer that $\cH_{\mathcal{G}}^*=\cH^*$.
Combining this with \eqref{4.29} and \eqref{4.30} and taking $\ve=1$, we obtain \eqref{4.27} and accomplish the proof.
\eo

Finally, Lemma \ref{le4.8} guarantees that when $|b|\leq\hat b$ and $\cM(g)\leq\hat M$, both $\sK_0$ and $\sK_1$ are nonempty.
By Theorem \ref{th3.3}, we have at least two recurrent solutions $\gam_0$ and $\gam_1$ for the problem \eqref{1.1} -- \eqref{1.3}, such that
$$\{(\gam_i(t),\theta_tg):t\in\R\}\ss\sK_i,\Hs i=0,\,1.$$
We finish the proof of Theorem \ref{th4.1} now.

\br Theorem \ref{th3.3} can be successfully applied to many other equations and more recurrent solutions may be obtained.

As an example, we consider the Dirichlet boundary-value problem
\be\label{4.31}u_t-\De u-au+u^3=g(t,x),\mb{ in }\W\X(\tau,\8),\hs u|_{\pa\W}=0\mb{ in }[\tau,\8),\ee
where $\W\subset\R$ is an open interval.
When $a>\mu_1$, we know the elliptic equation $-\De u-au+u^3=0$ with Dirichlet boundary value possesses at least three equilibria (see \cite[Theorem 2.44]{Car}), saying $0$ and $\pm u_0$, where $u_0$ is a positive solution.
One can easily see that $u_0$ (or $-u_0$) is a local attractor of the evolutionary problem \eqref{4.31} with $g\equiv0$.
Thus by similar discussion as above, the problem \eqref{4.31} with \eqref{1.3} has at least three recurrent solutions with distinct initial values, when $a>\mu_1$ and $g$ is a recurrent function with $\sup_{t\in\R}\|g(t,\cdot)\|$ sufficiently small.\er

\br Essentially, recurrent solutions of nonautonomous differential equation can be thought as solutions bifurcating from compact invariant sets of the corresponding autonomous equation.
In particular, if the autonomous equation has an evident Morse decomposition, we can obtain a lower bound of the number of recurrent solutions for the nonautonomous equation when the norm supremum of the recurrent forcing term is sufficiently small.
However, one necessarily notice that, only compact invariant sets whose Conley indices are not $\ol0$ count for calculating this number.
\er
%%%%%%%%%%%%%%%%%%%%%%%%%%%%%%%%%%%%%%%%%%%%%%%%%%%%%%%%%%


\begin{thebibliography}{99}
\small
\bibitem{B}G. D. Birkhoff, {\it Dynamical Systems}, Amer. Math. Soc. Colloq. Publ., vol. 9, Amer. Math. Soc., Providence, RI, 1966.
\bibitem{BCD}V. P. Bongolan-Walsh, D. Cheban, and J. Duan, Recurrent Motions in the Nonautonomous Navier-Stokes System, {\it Discret. Contin. Dyn. Syst-Ser. B}, 3: 255 - 262 (2003).
\bibitem{Car} A. N. Carvalho, J. A. Langa, and J. C. Robinson, {\it Attractors of Infinite Dimensional Nonautonomous Dynamical Systems}, Applied Mathematical Sciences Vol. 182, Springer New York, etc., 2013.
\bibitem{Ch}D. N. Cheban, {\it Global Attractors of Nonautonomous Dynamical Systems (Russian)}, Kishinev, State University of Moldova, 2002.
\bibitem{Cho1}Y. Choi, Dynamical bifurcation of the one dimensional modified Swift-Hohenberg equation, {\it Bull. Korean. Math. Soc.}, 52(4): 1241 - 1252 (2015).
\bibitem{Cho2}Y. Choi, T. Ha, J. Han, D. S. Lee, Bifurcation and final patters of a modified Swift-Hohenberg equation, {\it Discret. Contin. Dyn. Syst-Ser. B}, 22(7): 2543 - 2567 (2017).
\bibitem{Con1}C. Conley, {\it Isolated invariant sets and the Morse index}, Regional Conference Series in Mathematics, vol. 38, Amer. Math. Soc., Providence, RI, 1978.
\bibitem{Duan}N. Duan, W. Gao, Optimal control of a modified Swift-Hohenberg equation, {\it Electron. J. Differ. Equ.}, 2012(155): 1 - 12 (2012.)
\bibitem{GuoGuo1}C. Guo, Y. Guo, C. Li, Dynamic behaviors of a local modified stochastic Swift-Hohenberg equation with multiplicative noise, {\it Bound. Value Probl.} 2017(1): 1 - 13 (2017).
\bibitem{Hen}D. Henry, {\it Geometric theory of semilinear parabolic equations}, Lect. Notes in Math. 840, Springer Verlag, Berlin New York, 1981.
\bibitem{Kloeden} P. E. Kloeden, M. Rasmussen, {\it Nonautonomous Dynamical Systems}, Mathematical Surveys and Monographs 176, Amer. Math. Soc., 2011.
\bibitem{LaQ}R. E. LaQuey, S. M. Mahajan, P. H. Rutherford, et al, Nonlinear saturation of the trapped-ion mode, {\it Phys. Rev. Lett.}, 34(7): 391 - 394 (1975).
\bibitem{LM}J. Lega, J. V. Moloney, A. C. Newell, Swift-Hohenberg equation for lasers, {\it Phys. Rev. Lett.}, 73(22): 2978 - 2981 (1994).
\bibitem{LiW}D. S. Li, J. Y. Wei, J. T. Wang, On the dynamics of abstract retarded evolution equations, {\it Abstract Appl. Anal.}, 2013 (2013).
\bibitem{Park}S. H. Park, J. Y. Park, Pullback attractor for a non-autonomous modified Swift-Hohenberg equation, {\it Comput. Math. Appl.} 67(3): 542 - 548 (2014).
\bibitem{Pol}M. Polat, Global attractor for a modified Swift-Hohenberg equation, {\it Comput. Math. Appl.}, 57(1): 62 - 66 (2009).
\bibitem{PoM}Y. Pomeau and P. Manneville, Wave length selection in cellular flows, {\it Phys. Lett. A}, 75(4): 296 - 298 (1980).
\bibitem{Ryb}K. P. Rybakowski, {\it The homotopy index and partial differential equations}. Springer-Verlag, etc., 1980.
\bibitem{Se1}G. R. Sell, {\it Topological Dynamics and Ordinary Differential Equations}. Van Nostrand-Reinbold, London, 1971.
\bibitem{Se2}G. R. Sell and Y. You, {\it Dynamics of Evolutionary Equations}. Springer, 2002.
\bibitem{Sh2}B. A. Shcherbakov, Recurrent solutions of differential equations and general theory of dynamical systems, {\it Differential Equations (English Transl.)}, 3: 758 - 763 (1967).
\bibitem{Siva} G. L. Sivashinsky, Nonlinear analysis of hydrodynamic instability in laminar flames --- I. Derivation of basic equations, {\it Acta Astron.}, 4(11-12): 1177 - 1206 (1977).
\bibitem{Song}L. Song, Y. Zhang, T. Ma, Global attractor of a modified Swift-Hohenberg equation in $H^k$ spaces, {\it Nonlinear Anal.-Theory Methods Appl.}, 72(1): 183-191 (2010).
\bibitem{Sun}B. Sun, Optimal disturbed contral problem for the modified Swift-Hohenberg equations, {Electron. J. Differ. Equ.}, 2018(131): 1 - 13 (2018).
\bibitem{Swi}J. B. Swift, P. C. Hohenberg, Hydrodynamic fluctuations at the convective instability, {\it Phys. Rev. A}, 15(1): 319 - 328 (1977).
\bibitem{Tem}R. Temam, {\it Infinite-Dimensional Dynamical Systems in Mechanics and Physics} (Second Edition). Springer-Verlag, New York, 1997.
\bibitem{Wang}Z. Wang, X. Du, Pullback attractors for modified Swift-Hohenberg equation on unbounded domains with non-autonomous deterministic and stochastic forcing terms, {\it J. Appl. Anal. Comput.}, 7(1): 207 - 223 (2017).
\bibitem{Xiao}Q. Xiao, H. Gao, Bifurcation analysis of a modified Swift-Hohenberg equation, {\it Nonlinear Anal.-Real World Appl.}, 11: 4451 - 4464 (2010).
\bibitem{XuMa}L. Xu, Q. Ma, Existence of the uniform attractors for a non-autonomous modified Swift-Hohenberg equation, {\it Adv. Differ. Equ.}, 2015(1): 1 - 11 (2015).
\bibitem{Zhao}X. Zhao, Z. H. Ning, Optimal control problem for the modified Swift-Hohenberg equation in 3D case, {\it Math. Meth. Appl. Sci.}, 38(18): 4650 - 4662 (2015).
\bibitem{Zheng}J. Zheng, Optimal controls of multidimensional modified Swift-Hohenberg equation, {\it Int. J. Control}, 88(10): 1 - 9 (2015).
\end{thebibliography}
\end{document}